\documentclass[11pt]{amsart}
\usepackage{mathrsfs,amsthm,amssymb,amsmath,url}

\theoremstyle{plain}
    \newtheorem{thm}{Theorem}
    
    \newtheorem{lem}[thm]{Lemma}
    \newtheorem{prop}[thm]{Proposition}
    \newtheorem{cor}[thm]{Corollary}
    
\theoremstyle{definition}
    \newtheorem{defn}[thm]{Definition}
    \newtheorem{nota}[thm]{Notation}
\theoremstyle{remark}
    
\theoremstyle{remark}

\newcommand{\To}{\rightarrow}

\newcommand{\mult}{\times}
\newcommand{\ow}{\text{otherwise}}
\newcommand{\nin}{\notin}
\renewcommand{\ss}{\subseteq}
\newcommand{\sm}{\setminus}
\renewcommand{\c}{X \setminus}

\newcommand{\nf}{_{n_f}}

\newcommand{\cl}[1]{\langle #1 \rangle}

\newcommand{\OO}{{\mathscr O}}
\renewcommand{\O}{{\mathscr O}}
\newcommand{\On}{{\mathscr O}^{(n)}}

\newcommand{\Oo}{{\mathscr O}^{(1)}}

\DeclareMathOperator{\pol}{Pol}
\newcommand{\p}[1]{\pol(#1)}
\newcommand{\inv}{^{-1}}

\newcommand{\C}{{\mathscr C}}

\newcommand{\F}{{\mathscr F}}
\newcommand{\J}{{\mathscr J}}
\newcommand{\I}{{\mathscr I}}
\newcommand{\A}{{\mathscr A}}
\newcommand{\B}{{\mathscr B}}

\newcommand{\E}{{\mathscr E}}
\newcommand{\G}{{\mathscr G}}
\renewcommand{\S}{{\mathscr S}}
\newcommand{\M}{{\mathscr M}}

\newcommand{\U}{{\mathscr U}}
\renewcommand{\H}{{\mathscr H}}
\renewcommand{\L}{{\mathscr L}}
\newcommand{\X}{{\mathscr X}}

\newcommand{\un}{^{(n)}}

\newcommand{\uo}{^{(1)}}
\newcommand{\ut}{^{(2)}}

\newcommand{\rest}{\upharpoonright}

\DeclareMathOperator{\Const}{Const}

\author[M.\,Pinsker]{Michael Pinsker}
\address{Algebra\\TU Wien\\Wiedner Hauptstrasse 8-10/104\\A-1040 Wien, Austria}
\email{marula@gmx.at} \urladdr{http://www.dmg.tuwien.ac.at}
\title[Maximal clones containing the permutations]{Maximal clones on uncountable sets that include all permutations}


\thanks{The author would like to thank L. Heindorf for his many comments on earlier versions of the paper.
Support by DOC [Doctoral Scholarship Programme of the Austrian
Academy of Sciences] is gratefully acknowledged.}

\begin{document}

\begin{abstract}
    We first
    determine the maximal clones on a set $X$ of infinite regular cardinality $\kappa$ which contain all
    permutations but not all unary functions, extending a result of Heindorf's for countably
    infinite $X$. If $\kappa$ is countably infinite or weakly compact, this yields a list of
    all maximal clones containing the permutations since in that case the maximal
    clones above the unary functions are known. We then generalize a result
    of Gavrilov's to obtain on all infinite $X$ a list of all maximal submonoids of the monoid
    of unary functions which contain the permutations.
\end{abstract}
\maketitle\thispagestyle{empty}

\begin{section}{Clones and the Results}

\subsection{The clone lattice}

Let $X$ be a set of size $|X|=\kappa$. For each natural number
$n\geq 1$ we denote the set of functions on $X$ of arity $n$ by
$\On$. We set $\O=\bigcup_{n=1}^\infty \On$ to be the set of all
finitary operations on $X$. A \emph{clone} is a subset of $\O$
which contains the projection maps and which is closed under
composition. Since arbitrary intersections of clones are obviously
again clones, the set of clones on $X$ forms a complete algebraic
lattice $Cl(X)$ which is a subset of the power set of $\O$. The
clone lattice is countably infinite if $X$ has exactly two
elements, and of size continuum if $X$ is finite and has at least
three elements. On infinite $X$ we have $|Cl(X)|=2^{2^\kappa}$.

The dual atoms of the clone lattice are called \emph{maximal}
clones. On finite $X$ there exist finitely many maximal clones and
an explicit list of those clones has been provided by Rosenberg
\cite{Ros70}. Moreover, the clone lattice is dually atomic in that
case, that is, every clone is contained in a maximal one. For $X$
infinite the number of maximal clones equals the size of the whole
clone lattice (\cite{Ros76}, see also \cite{GS02}), so that it
seems impossible to know all of them. It has also been shown
\cite{GS04} that if the continuum hypothesis holds, then not every
clone on a countably infinite set is contained in a maximal one.

\subsection{Clones containing the bijections}

    However, even on infinite $X$ the sublattice of $Cl(X)$ of clones containing the set $\S$ of all permutations of
    $X$ is dually atomic since $\OO$ is finitely generated over
    $\S$: Call a set $A\ss X$ \emph{large} iff $|A|=|X|=\kappa$ and \emph{small} otherwise.
    Moreover, $A$ is \emph{co-large} iff $X\sm A$ is large, and
    \emph{co-small} iff $X\sm A$ is small. Set
    $$
        \I=\{f\in\Oo: f \text{ is injective and } f[X]
        \text{ is co-large}\}
    $$
    and
    $$
        \J=\{g\in\Oo: g^{-1}[y] \text{ is large for all } y\in
        X\}.
    $$
    It is readily verified that for arbitrary fixed $f\in\I$ and
    $g\in\J$ we have
    $$
        \I=\{\alpha\circ f: \alpha\in\S\}\text{ and
        } \J=\{\alpha\circ g\circ\beta: \alpha,\beta\in\S\}.
    $$
    Moreover,
    $$
        \Oo=\{j\circ i: j\in \J,i\in\I\}.
    $$
    Together with the well-known fact that $\Oo\cup\{p\}$ generates $\O$ for any
    binary injection $p$ we conclude that $\O$ is generated by $\S\cup\{p,f,g\}$.
    Hence Zorn's lemma implies that the interval $[\S,\O]$ is
    dually atomic.

    We will determine all maximal
    clones $\C$ on a base set of regular cardinality for which $\S\ss\C$ but not $\Oo\ss\C$. This has
    already been done for countable base sets by Heindorf in the article
    \cite{Hei02} in the following way: Let $\rho\subseteq X^J$ be a relation on $X$ indexed by $J$ and let
    $f\in\On$. We say that $f$ \emph{preserves} $\rho$ iff for all $r^1=(r^1_i:i\in J),\ldots,r^n=(r^n_i:i\in J)$ in
    $\rho$ we have $(f(r^1_i,\ldots,r^n_{i}):i\in J) \in \rho$. We define the clone of \emph{polymorphisms} $\pol(\rho)$
    of $\rho\ss X^J$ to consist exactly  of the
    functions in $\OO$ preserving $\rho$. In particular, if
    $\rho\subseteq X^{X^k}$ is a set of $k$-ary functions, then the
    polymorphisms of $\rho$ are exactly those $f\in\On$ for which the composite $f(g_1,\ldots,g_n)\in\rho$
    whenever $g_1,\ldots,g_n\in\rho$. Now it is obvious that since
    clones are closed under composition we have
    $\C\ss\p{\C\un}$ for any clone $\C$ and for all $n\geq 1$, where $\C\un=\C\cap\On$. Moreover,
    $\p{\C\un}\un=\C\un$. Therefore, if $\C$ is a maximal clone such
    that $\S\ss\C\uo\subsetneqq\Oo$, then
    $\C\ss\p{\C\uo}\subsetneqq\O$ holds. Hence $\C=\p{\C\uo}$ by
    the maximality of $\C$. We conclude that all
    maximal clones with $\S\ss\C\uo\subsetneqq\Oo$ are of the form
    $\p{\G}$, where $\S\ss\G\subsetneqq\Oo$ is a \emph{submonoid} of
    $\Oo$, that is, a set of unary functions closed under composition and containing the identity map.

    We say that a property holds for \emph{almost all} $y\in X$ iff
    the set of all elements for which the property does not hold
    is small. For $\lambda\leq \kappa$ a cardinal define a
    unary function $f$ to be \emph{$\lambda$-surjective} iff $|X\sm f[X]|<\lambda$. Instead of $\kappa$-surjective we also
    say \emph{almost surjective}; this means that the range of $f$ is co-small. $f$ is
    \emph{$\lambda$-injective} iff
    $|\{x\in X:\exists y\neq x\, (f(x)=f(y))\}|<\lambda$. For $\lambda=1$ or infinite, this is
    the case iff there exists a set $A\ss X$ such that
    $|A|<\lambda$ and such that the restriction of $f$ to the
    complement of $A$ is injective. \emph{Almost injective} means
    $\kappa$-injective.

    We are going to prove the following theorem.
    \begin{thm}\label{THM:bijections}
        Let $X$ be a set of regular cardinality $\kappa$. The maximal clones
        over $X$ which contain all bijections but not all unary
        functions are exactly those of the form $\pol(\G)$, where
        $\G\in\{\A,\B,\E,\F\}\cup\{\G_\lambda: 1\leq\lambda\leq\kappa,\,\lambda\text{ a cardinal}\}$ is one of the following submonoids of $\Oo$:
        \begin{enumerate}
            \item{$\A=\{f\in\Oo:f^{-1}[y]$ is small for
            almost all $y\in X\}$}
            \item{$\B=\{f\in\Oo:f^{-1}[y]$ is small for all $y\in X\}$}
            \item{$\E=\{f\in\Oo:f$ is almost surjective$\}$}
            \item{$\F=\{f\in\Oo:f$ is almost surjective or constant$\}$}
            \item{$\G_\lambda=\{f\in\Oo:$ if $A\subseteq X$ has cardinality $\lambda$ then $|X\setminus f[X\setminus A]|\geq \lambda\}$}
        \end{enumerate}
    \end{thm}

    \begin{cor}\label{COR:numberOfMaximal}
        Let $X$ be a set of regular cardinality $\kappa=\aleph_\alpha$.
        Then there exist $\max(|\alpha|,\aleph_0)$ maximal clones on $X$
        which contain all bijections but not all unary functions.
    \end{cor}

    On some infinite sets, namely countably infinite ones and sets of weakly compact cardinality, it is known that
    there exist exactly two maximal clones $\pol(T_1)$ and $\pol(T_2)$ which contain $\Oo$ ($T_1$ and $T_2$
    are certain sets of binary functions). See Gavrilov \cite{Gav65} for
    the countable, and Goldstern and Shelah \cite{GS02} for the uncountable.
    Hence in those cases, our theorem completes the list of maximal clones above $\S$.
    It is a fact that weakly compact cardinals $\kappa$
    satisfy $\kappa=\aleph_\kappa$. Thus we have

    \begin{cor}
        Let $X$ be a set of countably infinite or weakly compact cardinality $\kappa$.
        Then there exist $\kappa$ maximal clones
        which contain all bijections.
    \end{cor}

    Unfortunately things are not always that easy, as demonstrated by Goldstern and Shelah in
    \cite{GS02}: For many regular cardinalities of $X$, in particular for all
    successors of uncountable regulars,
    there exist $2^{2^{\kappa}}$ maximal clones which contain $\Oo$. It is interesting that whereas above $\Oo$ the
    number of maximal clones varies heavily with the partition
    properties of the underlying base set (2 for weakly compact
    cardinals, $2^{2^\kappa}$ for many others), the number of
    maximal clones above the permutations but not above $\Oo$ is a
    monotone function of $\kappa$ and always relatively small
    ($\leq\kappa$).

\subsection{Maximal submonoids of $\Oo$}
    Not all monoids appearing in Theorem
    \ref{THM:bijections} are maximal submonoids of $\Oo$ (by a maximal submonoid of $\Oo$ we mean a dual atom
    in the lattice of submonoids of $\Oo$ with inclusion).
    More surprisingly, there exist maximal submonoids of $\Oo$ above the
    permutations whose polymorphism clone is not maximal. Observe that submonoids of $\Oo$
    differ only formally from \emph{unary clones}, that is clones consisting
    only of essentially unary functions, and that the lattice of
    monoids which contain the permutations is dually atomic by the
    argument we have seen before. We are going to prove the
    following theorem for arbitrary infinite sets in the last
    section, generalizing a theorem due to Gavrilov
    \cite{Gav65} for countable base sets.

    \begin{thm}\label{THM:allMaximalMonoids}
        Let $X$ be an infinite set. If $X$ has regular cardinality, then the
        maximal submonoids of $\Oo$ which contain the permutations
        are exactly the monoid $\A$ and the monoids $\G_\lambda$ and $\M_\lambda$ for
        $\lambda=1$ and $\aleph_0\leq \lambda\leq\kappa$, $\lambda$ a cardinal, where
        $$
            \M_\lambda=\{f\in\Oo: f \text{ is
            }\lambda\text{-surjective or not
            }\lambda\text{-injective}\}.
        $$
        If $X$ has singular cardinality, then the same is true
        with the monoid $\A$ replaced by
        $$
            \A'=\{f\in\Oo: \exists \lambda < \kappa \,\,(\,|f\inv[\{x\}]|\leq\lambda\text{ for almost all }x\in X\,)\,\}.
        $$
    \end{thm}
    \begin{cor}\label{COR:numberOfMaximalMonoids}
        On a set $X$ of infinite cardinality $\aleph_\alpha$ there
        exist $2\,|\alpha|+5$ maximal submonoids of $\Oo$ that contain
        the permutations. Hence the smallest cardinality on which there
        are infinitely many such monoids is $\aleph_{\omega}$.
    \end{cor}

    Observe that the statement about singular cardinals in Theorem \ref{THM:allMaximalMonoids} differs only
    slightly from the corresponding one for regulars.
    We do not know whether Theorem \ref{THM:bijections} can be generalized to singulars, but
    in our proof we do use the
    regularity condition (in Proposition \ref{PROP:PolAisMaximal},
    Lemma \ref{LEM:dasGammaLemma} and permanently in Section \ref{SEC:almostUnary}).

\subsection{An equivalent definition of $\G_\lambda$}
In the countable versions of Theorems \ref{THM:bijections} and
\ref{THM:allMaximalMonoids} a different but equivalent definition
of $\G_{\aleph_0}$ was used: Define for $\lambda=1$ and for all
$\aleph_0\leq\lambda \leq\kappa$ monoids
$$
    \delta(\lambda)=\{f\in\Oo: f \text{ is }\lambda\text{-injective or not }\lambda\text{-surjective}\}
$$
(this definition and notation is due to Rosenberg \cite{Ros74}).
Then we have
\begin{lem}\label{LEM:whereHasDgone}
    $\delta(\lambda)=\G_\lambda$ for $\lambda=1$ and
    $\aleph_0\leq\lambda\leq\kappa$.
\end{lem}
\begin{proof}
    Note that for $\lambda=1$, $\lambda$-injective simply means
    injective and $\lambda$-surjective means surjective. The lemma
    is easily verified for that case, and we prove it for
    $\lambda$ infinite.\\
    Assuming $f\in\delta(\lambda)$ we show $f\in\G_\lambda$. It is
    clear that if $f$ is not $\lambda$-surjective, then
    $f\in\G_\lambda$. So assume $f$ is $\lambda$-surjective; then
    by the definition of $\delta(\lambda)$, $f$ is
    $\lambda$-injective. Now let $A\ss X$ be an arbitrary set of
    size $\lambda$. Assume towards contradiction that $|X\sm f[X\sm
    A]|<\lambda$. Then two things can happen: If $|f[A]\cap f[X\sm
    A]|\geq\lambda$, then $|\{x\in X:\exists y\neq
    x\,(f(x)=f(y))\}|\geq |\{x\in A:\exists y\in X\sm A\,
    (f(x)=f(y))\}|\geq\lambda$, contradicting the
    $\lambda$-injectivity of $f$. Otherwise, $A$ is mapped onto a
    set of size smaller than $\lambda$, again in contradiction to $f$ being
    $\lambda$-injective.\\
    To see the other inclusion, take any $f\nin\delta(\lambda)$.
    Then $f$ is not $\lambda$-injective; thus we can find $A\ss
    X$ of size $\lambda$ such that $f[X]=f[X\sm A]$. But then
    $|X\sm f[X\sm A]|$=$|X\sm f[X]|<\lambda$ as $f$ is
    $\lambda$-surjective. Hence, $f\nin\G_\lambda$.
\end{proof}

 Before we start with the proofs we fix some global notation.
\subsection{Notation}
    For a set of functions $\F$ we shall denote the smallest
    clone containing $\F$ by $\langle \F \rangle$. We call the projections which every clone contains
    $\pi^n_i$, where $n\geq 1$ and $1\leq i\leq n$.
    We write $n_f$ for the arity of a function
    $f\in\O$ whenever that arity has not yet been given another name. If $a\in X^n$ is an
    $n$-tuple and $1\leq k\leq n$ we write $a_k$ for the $k$-th component of
    $a$. The image of a set $A\ss X^n$ under a function $f\in\On$ we
    denote by $f[A]$. Similarly we write $f\inv[A]$ for the
    preimage of $A\ss X$ under $f$. If $A=\{c\}$ is a singleton we
    cut short and write $f\inv[c]$ rather than $f\inv[\{c\}]$.
    Occasionally we shall denote the constant function with value
    $c\in X$ also by $c$.
    Whenever we identify $X$ with its cardinality we let $<$ and
    $\leq$ refer to the canonical well-order on $X$.
\end{section}

\begin{section}{The proof of Theorem \ref{THM:bijections}}

    In this section we are going to prove Theorem
    \ref{THM:bijections}; it will be the direct consequence of Propositions \ref{PROP:PolAisMaximal}, \ref{PROP:PolGUnderPolA},
    \ref{PROP:PolGUnderPolB}, \ref{PROP:someIlambdaIsBelow},
    \ref{PROP:PolGlambdaIsMax}, \ref{PROP:GunderE},
    \ref{PROP:GunderF}, and \ref{PROP:remainingPolsNotMaximal}.
    The first part of the proof (Section \ref{SEC:asInCountable}) is not much more than a translation of
    Heindorf's paper \cite{Hei02} to arbitrary regular cardinals; the
    reader familiar with that article should not be surprised to find
    the same constructions here. In Section \ref{SEC:almostUnary}
    we have to go an own way to finish the proof.

\subsection{The beginning of the proof}\label{SEC:asInCountable}

We start with a general observation which will be useful.
\begin{lem}\label{LEM:unaryFunctionsSuffice}
    Let $\G$ be a proper submonoid of $\Oo$ such that
    $\cl{\pol(\G)\cup\{h\}}=\O$ for all unary $h\nin\G$. Then
    $\pol(\G)$ is maximal.
\end{lem}
\begin{proof}
    Let $f\nin \pol(\G)$ be given. Then there exist
    $h_1,\ldots,h_{n_f}\in\G$ such that $h=f(h_1,\ldots,h_{n_f})\nin\G$.
    Now $h\in\cl{\G\cup\{f\}}\subseteq\cl{\pol(\G)\cup\{f\}}$ and
    $\cl{\pol(\G)\cup\{h\}}=\O$ by assumption so that we conclude
    $\cl{\pol(\G)\cup\{f\}}=\O$.
\end{proof}

\subsubsection{The monoids $\A$ and $\B$}
\begin{prop}\label{PROP:PolAisMaximal}
    The clones $\p{\A}$ and $\p{\B}$ are maximal.
\end{prop}
\begin{proof}
    The maximality of $\p{\A}$ has been proved in
    \cite{Gav65} for the countable case and in \cite{Ros74} (Proposition 4.1) for
    arbitrary sets of regular cardinality (although not stated there, the regularity condition is necessary, for
    otherwise $\A$ is not closed under composition).

    For the maximality of $\p{\B}$, let a unary $h\nin\B$ be given; by
    Lemma \ref{LEM:unaryFunctionsSuffice}, it suffices to show
    $\cl{\pol(\G)\cup\{h\}}=\O$. By the definition of $\B$ there
    exists $c\in X$ such that the preimage $Y=h\inv [c]$ is large.
    Choose any injection $g:X\To Y$; then $h\circ g(x)=c$ for all $x\in X$.

    Now let $f\in\On$ be an arbitrary function and consider
    $\tilde{f}\in\O^{(n+1)}$ defined by
    $$
        \tilde{f}(x_1,\ldots,x_n,y)=\begin{cases}f(x_1,\ldots,x_n)&,y=c\\y&,y\neq
        c.\\\end{cases}
    $$
    We claim that $\tilde{f}\in\p{\B}$. For let
    $\alpha_1,\ldots,\alpha_n,\beta\in\B$ and $d\in X$ be given. If
    $\tilde{f}(\alpha_1,\ldots,\alpha_n,\beta)(x)=d$, then by the definition of $\tilde{f}$ either
    $\beta(x)=c$ and $f(\alpha_1(x),\ldots,\alpha_n(x))=d$ or $\beta(x)\neq
    c$ and $\beta(x)=d$. But since $\beta\in\B$, the set of all
    $x\in X$ such that $\beta(x)=c$ or $\beta(x)=d$ is small.
    Hence $\tilde{f}(\alpha_1,\ldots,\alpha_n,\beta)\inv[d]$ is small and so
    $\tilde{f}(\alpha_1,\ldots,\alpha_n,\beta)\in\B$.

    Now to finish the proof it is enough to observe that
    $f(x_1,\ldots,x_n)=\tilde{f}(x_1,\ldots,x_n,c)=
    \tilde{f}(x_1,\ldots,x_n,h\circ g(x_1))\in\cl{\p{\B}\cup\{h\}}$.
\end{proof}
    We will prove now that $\B$ is the only proper submonoid of
    $\A$ whose $\pol$ is maximal.

\begin{lem}\label{LEM:dasGammaLemma}
    If $f\nin\p{\A}$, then there exist $\alpha_1,\ldots,\alpha\nf\in\Oo$ constant or injective such that
    $f(\alpha_1,\ldots,\alpha\nf)\nin\A$.
\end{lem}
\begin{proof}
    Since $f\nin\p{\A}$, there exist $\beta_1,\ldots,\beta\nf\in\A$
    such that $f(\beta_1,\ldots,\beta\nf)\nin\A$. We will use induction over
    $n_f$. If $n_f=1$, then $f\nin\p{\A}\uo=\A$ so that
    $f(\pi^1_1)=f \nin\A$ which proves the assertion for that
    case.
    Now assume the lemma holds for all functions of arity at most
    $n_f-1$. Define for $1\leq i\leq n_f$ sets $B_i=\{y\in
    X:\beta_i\inv[y] \text{ is large}\}$. By definition of $\A$, all
    $B_i$ are small. Set
    $$
        \Gamma=(\beta_1,\ldots,\beta\nf)[X]\setminus\prod_{1\leq
        i\leq n_f}B_i\subseteq X^{n_f}
    $$
    \textbf{Claim.} There exists a large set $D\subseteq X$ such that $f\inv [d]\cap \Gamma$ is large for all $d\in
    D$.\\
    To prove the claim, set $D=\{d\in X: f(\beta_1,\ldots,\beta\nf)\inv[d]\text{ large}\}\setminus
    f[\prod_{1\leq i\leq n_f} B_i]$. The set $D$ is large as $f(\beta_1,\ldots,\beta\nf)\nin\A$ and as $\prod_{1\leq i\leq n_f}
    B_i$ is small. Define $A_d=(f(\beta_1,\ldots,\beta\nf))\inv[d]$ for each $d\in
    D$. Then $(\beta_1,\ldots,\beta\nf)[A_d]\subseteq \Gamma$ is
    large for all $d\in D$. Indeed, assume to the contrary that there exists $d\in D$ such that
    $(\beta_1,\ldots,\beta\nf)[A_d]$ is small; then, since $|X|=\kappa$ is regular, there is an $x\in
    (\beta_1,\ldots,\beta\nf)[A_d]$ so that
    $(\beta_1,\ldots,\beta\nf)\inv[x]$ is large. But then we would
    have $x\in\prod_{1\leq i\leq n_f}B_i$, in contradiction to the
    assumption that $d\nin f[\prod_{1\leq i\leq n_f} B_i]$.
    This proves the claim since $f\inv [d]\cap \Gamma=(\beta_1,\ldots,\beta\nf)[A_d]$ is
    large for every $d\in D$.

    Defining hyperplanes $H^i_b=\{x\in X^{n_f}: x_i=b\}$ for all $1\leq i\leq n_f$ and all $b\in X$, we can
    write $\Gamma$ as follows:
    $$
        \Gamma=(\bigcup_{i=1}^{n_f}\bigcup_{b\in
        B_i}\Gamma\cap H^i_b)\cup \Delta,
    $$
    where $\Delta=\Gamma\setminus \bigcup_{i=1}^{n_f}\bigcup_{b\in
        B_i} H^i_b$. Since $\kappa$ is regular and the union consists only of a small number of sets, we have that either there
    exist $1\leq i\leq n_f$ and some $b\in B_i$ such that $f\inv[d]\cap\Gamma\cap H^i_b$ is large for a large set
    of $d\in D$, or
    $f\inv[d]\cap\Delta$ is large for a large set of $d\in D$. We distinguish
        the two cases:\\
    \textbf{Case 1.} There exist $1\leq i\leq n_f$ and $b\in B_i$
    such that $f\inv [d]\cap\Gamma\cap H^i_b$ is large for many $d\in D$; say
    without loss of generality $i=n_f$. Then
    $f(\beta_1,\ldots,\beta_{n_f-1},b)\nin \A$. By induction
    hypothesis, there exist $\alpha_1,\ldots,\alpha_{n_f-1}$ injective
    or constant such that $f(\alpha_1,\ldots,\alpha_{n_f-1},b)\nin\A$.
    Setting $\alpha\nf(x)=b$ for all $x\in X$ proves the lemma.\\
    \textbf{Case 2.} $f\inv [d]\cap\Delta$ is large for many $d\in D$.
    Observe that for all $a\in X$ and all $1\leq i\leq n_f$, $\Delta\cap
        H^i_a$ is small, for otherwise $\beta_i\inv[a]$ would be
        large and thus $a\in B_i$, contradiction. Set
    $$
        C=\{c\in X: f\inv[c]\cap\Delta\text{ large}\}.
    $$
    By the assumption for this case, $C$ is large. Now fix any $g:
    X\To C$ such that $g\inv[c]$ is large for all $c\in C$. We
    define a function $\alpha:X\To\Delta$ such that
    $f\circ\alpha=g$; moreover, $\alpha_i=\pi^{n_f}_i\circ\alpha$ will be injective, $1\leq i\leq n_f$.
    Identify $X$ with its cardinality
    $\kappa$. Then all $\alpha_i$ are injective iff
    $\alpha_i(x)\neq\alpha_i(y)$ for all $y< x$ and all $1\leq i \leq
    n_f$. This is the case iff
    $$
        (\alpha_1,\ldots,\alpha\nf)(x)\in\Delta\setminus\bigcup_{y<x}\bigcup_{i=1}^{n_f}
        H^i_{\alpha_i(y)}.
    $$
    Using transfinite induction on $\kappa$, we define $(\alpha_1,\ldots,\alpha\nf)$ by picking
    $$
        (\alpha_1,\ldots,\alpha\nf)(x)\in
        (f\inv[g(x)]\cap\Delta)\setminus\bigcup_{y<x}\bigcup_{i=1}^{n_f}
        H^i_{\alpha_i(y)}.
    $$
    This is possible as $f\inv[g(x)]\cap\Delta$ is large for all
    $x\in X$
    whereas $\Delta\cap\bigcup_{y<x}\bigcup_{i=1}^{n_f}
    H^i_{\alpha_i(y)}$ is small. Clearly
    $f(\alpha_1,\ldots,\alpha\nf)=g\nin \A$ and the proof of the
    lemma is complete.
\end{proof}
\begin{prop}\label{PROP:PolGUnderPolA}
    Let $\G\subseteq \A$ be a submonoid of $\Oo$ which contains
    all permutations. Then either $\G\subseteq\B$ or
    $\p{\G}\subseteq\p{\A}$.
\end{prop}
\begin{proof}
    Assume $\G \nsubseteq \B$; we show $\p{\G}\subseteq\p{\A}$. Observe first that
    for all co-large $A\subseteq X$ and all $a\in
    X$ there exists $g\in\G$ such that $g[A]=\{a\}$.
    Indeed, choose any $h\in\G\setminus\B$. There exists $y\in X$
    such that $h\inv[y]$ is large. Choose bijections
    $\alpha,\beta\in\S$ with the property that $\alpha[A]\subseteq
    h\inv[y]$ and that $\beta(y)=a$. Then $g=\beta\circ
    h\circ\alpha$ has the desired property.\\
    Now let $f\nin\p{\A}$ be arbitrary; we show $f\nin\p{\G}$. By
    the preceding lemma there exist $\alpha_1,\ldots,\alpha_{n_f}$
    constant or injective such that
    $f(\alpha_1,\ldots,\alpha\nf)\nin\A$. Choose a large and co-large $A\subseteq
    X$ such that $f(\alpha_1,\ldots,\alpha\nf)\inv[x]\cap A$  is large for
    a large set of
    $x\in X$. We modify the $\alpha_i$ to $\gamma_i\in\G$ in
    such a way that $\alpha_i\rest_A=\gamma_i\rest_A$ for $1\leq i\leq n_f$: If $\alpha_i$ is injective, then we
    can choose $\gamma_i$ to be a bijection. If $\alpha_i$ is
    constant, then $\gamma_i$ is delivered by the observation we just made. Thus, as
    $f(\alpha_1,\ldots,\alpha\nf)\rest_A=f(\gamma_1,\ldots,\gamma\nf)\rest_A$ we have
    $f(\gamma_1,\ldots,\gamma\nf)\nin\A\supseteq\G$.
\end{proof}
\begin{prop}\label{PROP:PolGUnderPolB}
    Let $\G\subseteq \B$ be a submonoid of $\Oo$ which contains
    all permutations. Then $\p{\G}\subseteq\p{\B}$.
\end{prop}
\begin{proof}
    For arbitrary $f\nin\p{\B}$ we show $f\nin\p{\G}$. There
    are $\beta_1,\ldots,\beta\nf\in\B$ such that there exists $c\in
    X$ with the property that $f(\beta_1,\ldots,\beta\nf)\inv[c]$ is
    large. Define $\Gamma=(\beta_1,\ldots,\beta\nf)[X]$. Then since $\beta_i\in\B$,
    $H^i_a\cap\Gamma$ is small for all $1\leq i\leq n_f$ and all
    $a\in X$, where $H^i_a=\{x\in X^{n_f}:x_i=a\}$.
    Moreover, $f\inv[c]\cap\Gamma$ is large. Just like
    at the end of the proof of Lemma
    \ref{LEM:dasGammaLemma}, we can construct injective
    $\alpha_1,\ldots,\alpha\nf$ such that
    $f(\alpha_1,\ldots,\alpha\nf)$ is constant with value $c$. Choose
    $A\subseteq X$ large and co-large and bijections
    $\gamma_1,\ldots,\gamma\nf$ such that $\gamma_i\rest_A=\alpha_i\rest_A$ for $1\leq i\leq n_f$.
    Then, being constant on $A$,
    $f(\gamma_1,\ldots,\gamma\nf)\nin\B\supseteq\G$. Thus,
    $f\nin\p{\G}$.
\end{proof}
\subsubsection{Generous functions}
We now turn to monoids $\G\supseteq\S$ which are not submonoids of
$\A$. Our first goal is Proposition \ref{PROP:someIlambdaIsBelow},
in which we give a positive description of such monoids.

\begin{defn}
    A function $f\in\Oo$ is called \textit{generous} iff
    $f\inv[y]$ is either large or empty for  all $y\in X$.
\end{defn}
\begin{nota}
    Let $0\leq\lambda\leq \kappa$ be a cardinal. We denote by
    $\I_\lambda$ the set of all generous functions $f$ with the
    property that $|X\setminus f[X]|=\lambda$.
\end{nota}

The verification of the following simple facts is left to the
reader.

\begin{lem}\label{LEM:basicIlambdaFacts}
    \begin{enumerate}
        \item{If $g\in \Oo$ is generous, then $f\circ g$ is generous for all $f\in\Oo$.}
        \item{$\I_\lambda$ is a subsemigroup and $\I_\lambda\cup\S$ a submonoid of $\Oo$ for all $\lambda\leq\kappa$.}
        \item{If $\lambda<\kappa$ and $f,g\in \I_\lambda$, then there exist $\alpha,\beta\in\S$ such that $f=\alpha\circ g\circ\beta$.}
        \item{$\I_\kappa$ contains all generous functions with small range, in particular the constant functions.}
        \item{If $g\in \I_\kappa$ has large range, then $\cl{\S\cup\{g\}}\supseteq \I_\kappa$.}
    \end{enumerate}
\end{lem}

\begin{lem}\label{LEM:generousInfiniteRange}
    If $g\nin\A$, then there exists $\alpha\in\S$ such that the
    function $g\circ\alpha\circ g$ is generous and has large
    range.
\end{lem}
\begin{proof}
    There exists a large set $A\ss X$ such that $g\inv[a]$ is large for
    all $a\in A$. Fix any $a_0\in A$ and set $B_0=g\inv[a_0]$. Write $A$ as a disjoint
    union: $A=\{a_0\}\cup A_1\cup A_2$,
    with $A_1,A_2$ large, and set $B_i=g\inv [A_i]$, $i=1,2$. Choose any injective partial
    mapping $\widetilde{\alpha}$ defined on $X\sm A_2$ which
    satisfies $\widetilde{\alpha}[A_1]=B_1$ and
    $\widetilde{\alpha}[X\sm(A_1\cup A_2)]\subseteq B_0$. Since
    both domain and range of
    $\widetilde{\alpha}$ are large and co-large, we can extend the
    function to a bijection $\alpha\in\S$. Now $g\circ\alpha\circ
    g[X]\supseteq g\circ\alpha[A_1]=g[B_1]=A_1$, so
    $g\circ\alpha\circ g$ has large range. For all $x\in B_1\cup
    B_2$, the equivalence class of $x$ in the kernel of $g$ is
    large, and hence also the class of $x$ in the kernel of
    $g\circ\alpha\circ g$. If on the other hand $x\nin B_1\cup
    B_2$, then $g(x)\nin A_1\cup A_2$, so that $\alpha\circ g(x)\in
    B_0$ and therefore $g\circ\alpha\circ g(x)=a_0$. But for all
    $y\in B_0$ we have $g\circ\alpha\circ g(y)=a_0$ so that the
    kernel class of $x$ is again large. Whence, $g\circ\alpha\circ
    g$ is generous.
\end{proof}

\begin{prop}\label{PROP:someIlambdaIsBelow}
    Let $\G\subseteq\Oo$ be a monoid containing all bijections.
    Then either $\G\subseteq\A$ or there exists a cardinal
    $\lambda\leq\kappa$ such that $\I_\lambda\subseteq\G$.
\end{prop}
\begin{proof}
    This is an immediate consequence of Lemmas
    \ref{LEM:basicIlambdaFacts} and
    \ref{LEM:generousInfiniteRange}.
\end{proof}
The preceding proposition implies that when considering submonoids
$\G$ of $\Oo$ which contain the permutations, we can from now on
assume that $\I_\lambda\subseteq\G$ for some $\lambda$, since we
already treated the case $\G\subseteq\A$. We distinguish two cases
corresponding to the minimal $\lambda$ with the property that
$\I_\lambda\subseteq\G$, $\lambda>0$ and $\lambda=0$.

\subsubsection{The case $0<\lambda\leq\kappa$}

We shall now investigate the case where $\G\nsupseteq\I_0$ but
$\G$ contains $\I_\lambda$ for some $0<\lambda<\kappa$. The
following facts about the $\G_\lambda$ are left to the reader. The
proof of (4) can be found in \cite{Ros74} (Lemma 5.2).

\begin{lem}\label{LEM:GlambdaFacts} The following statements hold for all
$1\leq\lambda\leq\kappa$.
    \begin{enumerate}
    \item{If $g\in\On$ and $|X\sm g[X^n]|\geq\lambda$, then
    $g\in\p{\G_\lambda}$.}
    \item{$\G_\lambda$ is a submonoid of $\Oo$.}
    \item{$\G_n\supsetneqq\G_{n+1}$ for all $1\leq n< \aleph_0$.}
    \item{For $\lambda=1$ and for $\lambda\geq\aleph_0$,
    $\G_\lambda$ is a maximal submonoid of $\Oo$.}
    \end{enumerate}
\end{lem}

\begin{lem}\label{LEM:existsLambda0}
    Let $1\leq\lambda\leq\kappa$. If $h\nin\G_\lambda$, then there exists a $\lambda_0 <\lambda$
    such that $\cl{\I_\lambda\cup\S\cup\{h\}}\supseteq \I_{\lambda_0}$. In particular,
    $\cl{\G_\lambda\cup\{h\}}\supseteq \I_{\lambda_0}$.
\end{lem}
\begin{proof}
    There exists $A\subseteq X$, $|A|=\lambda$ such that $|X\sm
    h[X\sm A]|<\lambda$. Set $\lambda_0=|X\sm
    h[X\sm A]|$. Choose a generous function $g$ with $g[X]=X\sm
    A$. Then $g\in\I_\lambda$ since $|X\sm g[X]|=|A|=\lambda$;
    thus, $h\circ g\in\cl{\I_\lambda\cup\{h\}}$. On the other
    hand, $h\circ g\in \I_{\lambda_0}$ and hence $\cl{\I_\lambda\cup\S\cup\{h\}}\supseteq
    \I_{\lambda_0}$ by Lemma \ref{LEM:basicIlambdaFacts} (3). The
    second statement is a direct consequence of the inclusion
    $\G_\lambda\supseteq\I_\lambda\cup\S$.
\end{proof}
\begin{lem}\label{LEM:theStrangeG}
    Let $B\ss X$, $|B|=\lambda_0<\lambda\leq\kappa$, and let
    $g\in\O\ut$ such that $g$ maps $(X\sm B)^2$ bijectively onto
    $X$ and such that $|g[B\mult X]\cup g[X\mult B]|<\kappa$. Then
    $g\in\p{\G_\lambda}$.
\end{lem}
\begin{proof}
    Let $\alpha,\beta\in\G_\lambda$ be given, and take an
    arbitrary $A\subseteq X$ of size $\lambda$. We have to show
    $|X\sm g(\alpha,\beta)[X\sm A]|\geq\lambda$. For
    $C=X\sm\alpha[X\sm A]$ we have $|C|\geq\lambda$. Thus, there
    exists some $c\in C\sm B$. Obviously, $g(\alpha,\beta)[\c A]\ss
    g[(\c \{c\})\mult X]$. But the conditions on $g$ yield that $g[(\c \{c\})\mult X]$ and
    $g[\{c\}\mult (\c B)]\sm (g[X\mult B]\cup g[B\mult X])$ are
    disjoint. Since $|g[\{c\}\mult (\c B)]|=\kappa$ and $|g[X\mult B]\cup g[B\mult
    X]|<\kappa$, this implies that $g(\alpha,\beta)$ misses
    $\kappa$ values on $X\sm A$ and hence,
    $g(\alpha,\beta)\in\G_\lambda$ and $g\in\p{\G_\lambda}$.
\end{proof}
\begin{prop}\label{PROP:PolGlambdaIsMax}
    \begin{enumerate}
    \item{$\p{\G_\lambda}$ is a maximal clone for all
    $1\leq\lambda\leq\kappa$.}
    \item{Let $\G\ss\Oo$ be a monoid containing all bijections as well as some $\I_\lambda$,
    where $0\leq \lambda\leq\kappa$, and let $\lambda$ be minimal
    with this property. If $\lambda>0$, then
    $\p{\G}\ss\p{\G_\lambda}$.}
    \end{enumerate}
\end{prop}
\begin{proof}
    (1) We show $\cl{\p{\G_\lambda}\cup\{h\}}=\O$ for an arbitrary $h\in\Oo\sm\G_\lambda$.
    By Lemma \ref{LEM:existsLambda0}, there
    exists $\lambda_0<\lambda$ such that
    $\I_{\lambda_0}\ss\cl{\G_\lambda\cup\{h\}}$. Now choose $B$
    and $g\in\p{\G_\lambda}$ as in Lemma \ref{LEM:theStrangeG}.
    Consider $\alpha: X\To (X\sm B)^2$ such that $\alpha$ takes every value twice.
    Clearly, $\alpha_1=\pi^2_1\circ\alpha$ and $\alpha_2=\pi^2_2\circ\alpha$ are elements of $\I_{\lambda_0}$.
    The function $p=g(\alpha_1,\alpha_2)=g\circ\alpha$ maps $X$ onto $X$ and takes
    every value twice as well. Therefore we can find a co-large
    set $A$ such that $p[A]=X$. Now fix a mapping $q: X\To A$
    so that $p\circ q$ is the identity map on $X$. Let an
    arbitrary $f\in\O$ be given. Then $q\circ f[X^{n_f}]\subseteq
    A$ is co-large which immediately implies $q\circ
    f\in\p{\G_\lambda}$. But then $f=p\circ(q\circ f)=f\in\cl{\p{\G_\lambda}\cup\{h\}}$ and
    so $\cl{\p{\G_\lambda}\cup\{h\}}=\O$ as $f$ was arbitrary.\\
    (2) First we claim that $\G\ss\G_\lambda$. Indeed, assume there exists $h\in
    \G\sm\G_\lambda$. Then, as $\I_\lambda\cup\S\ss\G$, by Lemma \ref{LEM:existsLambda0}
    there exists $\lambda_0<\lambda$ such that
    $\I_{\lambda_0}\ss\G$, in contradiction to the minimality of
    $\lambda$.
    Now let $f\nin\p{\G_\lambda}$ be arbitrary; we prove
    $f\nin\p{\G}$. There exist
    $\alpha_1,\ldots,\alpha\nf\in\G_\lambda$ such that
    $f(\alpha_1,\ldots,\alpha\nf)\nin\G_\lambda$. That is, there
    exists $A\subseteq X$ of size $\lambda$ with the property that
    $|\c f[\Gamma]|<\lambda$, where
    $\Gamma=\{(\alpha_1(x),\ldots,\alpha\nf (x)):x\in\c A\}$. Since
    $\alpha_i\in\G_\lambda$, $1\leq i\leq n_f$, for each $i$ there
    exists a set $B_i\subseteq X$, $|B_i|=\lambda$, such that
    $\alpha_i[\c A]\cap B_i=\emptyset$. Then $\Gamma\subseteq
    \Delta=(\c B_1)\mult \ldots \mult (\c B\nf)$. Choose $\beta: X\To \Delta$ onto and generous.
    Clearly $\beta_i=\pi^{n_f}_i\circ\beta\in\I_\lambda\ss\G$
    for all $1\leq i\leq n_f$. Now we choose any $C\subseteq X$ of size
    $\lambda$ such that $\beta[X\setminus C]=\beta[X]$.
    This is possible since $\beta$ is generous. Then we have that $f(\beta_1,\ldots,\beta\nf)[\c
    C]=f[\Delta]\supseteq f[\Gamma]$ and so, as $|X\sm f[\Delta]|\leq|X\sm f[\Gamma]|<\lambda$,
    $f(\beta_1,\ldots,\beta\nf)\nin\G_\lambda\supseteq\G$. Hence,
    $f\nin\p{\G}$.
\end{proof}

\subsubsection{The case $\lambda=0$ and $\G\ss\F$}
In the following proposition we treat the case where
$\I_0\subseteq\G\ss\E\ss\F$.
\begin{prop}\label{PROP:GunderE}
    \begin{enumerate}
    \item{$\p{\E}$ is a maximal clone.}
    \item{If $\G\ss\Oo$ is a monoid containing all bijections as well as
    $\I_0$, and if $\G\ss\E$, then $\p{\G}\ss\p{\E}$.}
    \end{enumerate}
\end{prop}
\begin{proof}
    (1) We prove that for any unary $h\nin\E$ we have
    $\cl{\p{\E}\cup\{h\}}=\O$. By definition $h[X]$ is
    co-large, so we can fix $A\ss X$ large and co-large such that
    $A\cap h[X]=\emptyset$. Choose any $g\in\Oo$ which maps $A$
    onto $X$ and which is constantly $0\in X$ on $\c A$. Then
    $g\in\E$ as it is onto. Moreover, $g\circ h$ is constantly
    $0$. Now let an arbitrary $f\in\On$ be given and define a
    function $\tilde{f}\in\O^{(n+1)}$ by
    $$
        \tilde{f}(x_1,\ldots,x_n,y)=\begin{cases}f(x_1,\ldots,x_n)&,y=0\\y&,\ow\end{cases}
    $$
    Then $\tilde{f}\in\p{\E}$. Indeed, this follows from the
    inclusion
    $\tilde{f}(\alpha_1,\ldots,\alpha_n,\beta)[X]\supseteq\beta[X]\sm\{0\}$ for arbitrary
    $\alpha_1,\ldots,\alpha_n,\beta\in\Oo$. Now
    $f(x)=\tilde{f}(x,0)=\tilde{f}(x,g\circ h(x_1))$ for all
    $x\in X^n$ and so $f\in \cl{\p{\E}\cup\{h\}}$.\\
    (2) Taking an arbitrary $f\nin\p{\E}$ we show that
    $f\nin\p{\G}$. There exist $\alpha_1,\ldots,\alpha\nf$ almost
    surjective such that $f(\alpha_1,\ldots,\alpha\nf)$ is not almost
    surjective. Consider a small set $A\ss X$ so that
    $A\cup\alpha_i[X]=X$ for all $1\leq i\leq n_f$. Let $\gamma$
    be a surjection from $\c A$ onto $X$ and define for $1\leq i\leq
    n_f$ functions
    $$
        \beta_i(x)=\begin{cases}\alpha_i\circ\gamma(x)&,x\in\c
        A\\x&,x\in A\end{cases}
    $$
    Clearly, all $\beta_i$ are surjective and
    $f(\beta_1,\ldots,\beta\nf)[X]=f(\alpha_1,\ldots,\alpha\nf)[X]\cup\{f(x,\ldots,x):x\in
    A\}$ is co-large. Fix any $\delta\in\I_0$. Obviously
    $\beta_i\circ\delta\in\I_0\ss\G$ and also
    $f(\beta_1\circ\delta,\ldots,\beta\nf\circ\delta)[X]$ is
    co-large. Thus
    $f(\beta_1\circ\delta,\ldots,\beta\nf\circ\delta)\nin\E\supseteq\G$
    so that we infer $f\nin\p{\G}$.
\end{proof}

In a next step we see what happens in the case
$\I_0\subseteq\G\ss\F$ and $\G\nsubseteq\E$.

\begin{prop}\label{PROP:GunderF}
    \begin{enumerate}
    \item{$\p{\F}$ is a maximal clone.}
    \item{If $\G\ss\F$ is a monoid which contains $\I_0$ as well as all
    bijections, and if $\G\nsubseteq \E$, then
    $\p{\G}\ss\p{\F}$.}
    \end{enumerate}
\end{prop}

\begin{proof}
    (1) can be found in \cite{Ros74} (Proposition 3.1).\\
    For (2), let $f\nin\p{\F}$ and fix
    $\alpha_1,\ldots,\alpha\nf\in\F$ satisfying
    $f(\alpha_1,\ldots,\alpha\nf)\nin\F$. Since $\G\nsubseteq\E$ but $\G\ss\F$, $\G$ must contain a constant function,
    and hence all constant functions as $\S\ss\G$. For those of the $\alpha_i$
    which are not constant we construct $\beta_i$ as in
    the proof of the preceding proposition, and for the constant
    ones we set $\beta_i=\alpha_i$. Observe that it is impossible that all $\alpha_i$ are
    constant. Choosing any $\delta\in\I_0$
    we obtain that for all $1\leq i\leq n_f$, $\beta_i\circ\delta$ is either constant
    or an element of $\I_0$, and hence in either case an element
    of $\G$. But as in the preceding proof,
    $f(\beta_1\circ\delta,\ldots,\beta\nf\circ\delta)\nin\F\supseteq\G$
    so that $f\nin\p{\G}$.
\end{proof}

\subsubsection{The case $\lambda=0$ and $\G\nsubseteq\F$} To conclude, we
consider submonoids $\G$ of $\Oo$ which contain the bijections as
well as $\I_0$, but which are not submonoids of $\F$. It turns out
that the polymorphism clones of such monoids are never maximal. We
start with a simple fact about such monoids.

\begin{lem}
    Let $\G\ss\Oo$ be a monoid containing $\S\cup\I_0$ such that
    $\G\nsubseteq\F$. Then $\X=\{\rho\in\Oo:|\rho[X]|=2\text{ and }\rho \text{ is generous}\}\ss\G$.
\end{lem}
\begin{proof}
    Let $f\in\G\sm\F$. Since $f$ is not constant there exist
    $a\neq b$ in the range of $f$. Let $s:\c f[X]\To X$ be onto
    and generous and define $g\in\Oo$ by
    $$
        g(x)=\begin{cases}s(x)&,x\nin f[X]\\
                          a&,x=a\\
                          b&,\ow
        \end{cases}
    $$
    Then $g\in\I_0\ss\G$ and so $g\circ f\circ g\in \G$. On the other
    hand, $g\circ f\circ g\in\X$ which proves the lemma since obviously any function of
    $\X$ together with the permutations
    generate all of $\X$.
\end{proof}

\begin{nota}
    We set $\L=\cl{\X\cup\I_0\cup\S}$. Moreover, we write
    $\Const$ for the set of all constant functions.
\end{nota}
The following description of $\L$ is readily verified.
\begin{lem}
    $\L=\Const\cup\X\cup\I_0\cup\S$. In words, $\L$ consists
    exactly of the bijections as well as of all generous functions
    which are either onto or take at most two values.
\end{lem}

\begin{defn} A function $f(x_1,\ldots,x_n)\in\On$ is \emph{almost unary} iff
    there exist a mapping
    $F$ from $X$ to the power set of $X$ and some $1\leq k\leq n$ such that $F(x)$ is small for all
    $x\in X$ and such that for all
    $(x_1,\ldots,x_n)\in X^n$ we have $f(x_1,\ldots,x_n)\in F(x_k)$.
    We denote the set of all almost unary functions by
    $\U$.
\end{defn}

It is easy to see that on a base set of regular cardinality, $\U$
is a clone which contains $\Oo$. See \cite{Pin032} for a list of
all clones above $\U$; there are countably many, so in particular
$\U$ is not maximal. The reason for us to consider almost unary
functions is the following lemma.

\begin{lem}\label{LEM:LandNotAlmostUnary}
    Let $f\in\O\un\sm\U$ be any function which is not almost
    unary. Then $\cl{\{f\}\cup\L}\supseteq\Oo$.
\end{lem}

Therefore we have

\begin{prop}\label{PROP:remainingPolsNotMaximal}
    If $\G\ss\Oo$ is a nontrivial monoid such that $\S\cup\I_0\ss\G$ and
    such that $\G\nsubseteq\F$, then $\p{\G}\subseteq\U$. In particular, $\p{\G}$
    is not maximal.
\end{prop}

In his proof for countable base sets, Heindorf used the following
completeness criterion which has been shown by Gavrilov
\cite{Gav65} (Lemma 31 on page 51) to hold on countable base sets;
it will follow from our proof of Lemma
\ref{LEM:LandNotAlmostUnary} that this criterion holds on all
regular cardinals.

\begin{prop}\label{PROP:generalizedGavrilov}
    Let $X$ have regular cardinality. If $\G\ss\Oo$ is a monoid containing $\S\cup\I_0\cup\X$, and
    if $\H\ss\O$ is a set of functions such that
    $\cl{\Oo\cup\H}=\O$, then $\cl{\G\cup\H}=\O$.
\end{prop}

The criterion can be used alternatively to show that $\p{\G}$ is
not maximal for the remaining monoids $\G$: We have just seen that
$\X\ss\G$ so we can apply Proposition
\ref{PROP:generalizedGavrilov}. Suppose towards contradiction that
$\p{\G}$ is maximal. Since $\p{\G}\uo=\G\subsetneqq\Oo$ we
have $\cl{\Oo\cup\p{\G}}=\O$. But then setting $\H=\p{\G}$ in
the proposition yields that $\cl{\G\cup\p{\G}}=\O$, which is
impossible as $\cl{\G\cup\p{\G}}=\p{\G}\neq\O$, contradiction.

\subsection{The proof of Lemma \ref{LEM:LandNotAlmostUnary} and Proposition
\ref{PROP:generalizedGavrilov}.} \label{SEC:almostUnary}

\begin{lem}\label{LEM:universalFunction}
    Let $u\in\Oo$ be injective and not almost surjective. Then
    $\cl{\{u\}\cup\I_0} \supseteq \Oo$. In particular,
    $\cl{\{u\}\cup\L}\supseteq
    \Oo$.
\end{lem}
\begin{proof}
    Let an arbitrary $f\in\Oo$ be given. Take any $s: X\sm u[X] \To X$
    which is generous and onto. Now define $g\in \Oo$ by
    $$
        g(x)=\begin{cases}f(u\inv(x))&,x\in
        u[X]\\s(x)&,\ow\end{cases}
    $$
    Since $g\rest_{X\sm u[X]}=s$ we have $g\in\I_0$. Clearly, $f=g\circ
    u\in\cl{\{u\}\cup\I_0}$.
\end{proof}

Our strategy for proving Lemma \ref{LEM:LandNotAlmostUnary} is to
show that $\L$ together with a not almost unary $f$ generate a
function $u$ as in Lemma \ref{LEM:universalFunction}. We start by
observing that $\L$ and $f$ generate functions of arbitrary range.

\begin{lem}\label{LEM:largeAndColargeRange}
    Let $f\in\On\sm\U$. Then there exists a unary $g\in\cl{\{f\}\cup\L}$
    such that the
    range of $g$ is large and co-large.
\end{lem}
\begin{proof}
    We distinguish two cases.\\
    \textbf{Case 1.} For all $1\leq i\leq n$ and all $c\in X$ it
    is true that the image of the hyperplane $H^i_c$ under $f$ is
    co-small, where $H^i_c=\{x\in X^{n}: x_i=c\}$.
    Then consider an arbitrary large and co-large  $A\ss
    X$. Set
    $\Gamma=f\inv[X\sm A]\ss X^n$ and let $\alpha: X\To\Gamma$ be onto. By the assumption for this case,
    $f[H^i_c]\sm A$ is still
    large for all $1\leq i\leq n$ and all $c\in X$. Thus
    the components $\alpha_i=\pi^n_i\circ\alpha$ are generous and
    onto; hence, $\alpha_i\in\I_0\ss\L$ for all $1\leq i\leq n$. But now
    $f(\alpha_1,\ldots,\alpha_n)[X]=f[X^n]\sm A$ is large and
    co-large so that it suffices to set $g=f\circ\alpha$.\\
    \textbf{Case 2.} There exist $1\leq i\leq n$ and $c\in X$
    such that the image $f[H^i_c]$ of the hyperplane $H^i_c$ is co-large,
    say without loss of generality $i=1$. Since $f\nin\U$, there
    exists $d\in X$ satisfying that $f[H^1_d]$ is
    large. Choose $\Gamma\ss X^{n-1}$ large and co-large such that
    $f[\{d\}\mult\Gamma]$ is large and such that $f[H^1_c]\cup f[\{d\}\mult\Gamma]$ is still
    co-large. Take moreover $\alpha_2,\ldots,\alpha_n\in\I_0$ so that
    $(\alpha_2,\ldots,\alpha_n)[X]=X^{n-1}$. Now we define
    $\alpha_1\in\Oo$ by
    $$
        \alpha_1(x)=
        \begin{cases}d&,(\alpha_2,\ldots,\alpha_n)(x)\in\Gamma\\
        c&,\ow.
        \end{cases}
    $$
    Clearly, $\alpha_1\in\X\ss\L$. Now it is enough to set $g=f(\alpha_1,\ldots,\alpha_n)$ and observe
    that $g[X]=f[\{c\}\mult (X^{n-1}\sm\Gamma)]\cup f[\{d\}\mult\Gamma]$
    is large and co-large.
\end{proof}
\begin{lem}\label{LEM:arbitraryRange}
    Let $f\in\On\sm\U$. Then for all nonempty $A\subseteq X$ there
    exists $h\in\cl{\{f\}\cup\L}$ with $h[X]=A$.
\end{lem}
\begin{proof}
    By Lemma \ref{LEM:largeAndColargeRange} there exists $g\in \cl{\{f\}\cup\L}$
    having a large and co-large range. Now taking any
    $\delta\in\I_0\ss\L$ with $\delta[g[X]]=A$ and setting
    $h=\delta\circ g$ proves the assertion.
\end{proof}
\begin{lem}\label{LEM:allGenerousFunctions}
    If $f\in\On\sm\U$, then $\cl{\{f\}\cup\L}$ contains all generous
    functions.
\end{lem}
\begin{proof}
    Let any generous $g\in\Oo$ be given and take with the help of the preceding lemma $h\in\cl{\{f\}\cup\L}$ with
    $h[X]=g[X]$. By setting $h'=h\circ\delta$, where
    $\delta\in\I_0\ss\L$ is arbitrary, we obtain a generous
    function with the same property. Now it is clear that there
    exists a bijection $\sigma\in\S\ss\L$ such that
    $g=h'\circ\sigma$.
\end{proof}

Now that we know that we have all generous functions we want to
make them injective. We start by reducing the class of functions
$f$ under consideration.

\begin{lem}\label{LEM:columnsWithManyDifferentValues}
    If $f\in\On\sm\U$ is so that for all $1\leq i\leq n$ and for all $a,b\in X$
    the set of all tuples
    $(x_1,\ldots,x_{i-1},x_{i+1},\ldots,x_n)\in X^{n-1}$ with
    $f(x_1,\ldots,x_{i-1},a,x_{i+1},\ldots,x_n)\neq
    f(x_1,\ldots,x_{i-1},b,x_{i+1},\ldots,x_n)$ is small, then
    $\cl{\{f\}\cup\L}\supseteq\Oo$.
\end{lem}
\begin{proof}
    Since
    $f\nin\U$ we can for every $1\leq i\leq n$ choose $c_i\in X$
    such that $f[H^i_{c_i}]$ is large.
    Choose moreover for every $1\leq i\leq n$ large sets $A_i\ss f[H^i_{c_i}]$ such that
    $\bigcup_{i=1}^n A_i$ is co-large and such that $A_i\cap
    A_j=\emptyset$ for $i\neq j$. Write each $A_i$ as a disjoint
    union of many large sets: $A_i=\bigcup_{x\in X} A^x_i$. Let $\triangleleft$ be any well-order of
    $X^n$ of type $\kappa$. Define
    $\Gamma\ss X^n$ by $x\in\Gamma$ iff there exists $1\leq i\leq
    n$ such that $f(x)\in A_i^{x_i}$ and whenever $y\triangleleft
    x$ and $y\in\Gamma$ then $f(x)\neq f(y)$. Observe that the latter
    condition ensures that $f\rest_\Gamma$ is injective.

    Now observe that for all $1\leq i\leq n$, all $d\in X$ and all large $B\ss A_i$
    we have that $f[H^i_d]\cap B$ is
    large. Indeed, say without loss of generality $i=1$ and set
    $D=\{(x_2,\ldots,x_n):f(d,x_2,\ldots,x_n)\neq f(c_1,x_2,\ldots,x_n)\}$.
    Then $D$ is small by our assumption. Now $|f[H^1_d]\cap B|\geq |f[\{d\}\mult
    (X^{n-1}\sm D)]\cap B| = |f[\{c_1\}\mult
    (X^{n-1}\sm D)]\cap B|=\kappa$. In particular, this observation
    is true for $B=A_i^d$. This implies that the set
    $\{x\in\Gamma: x_i=d\}$ is large for all $1\leq i\leq n$ and
    all $d\in X$. Moreover, $\Gamma$ itself is large.

    Therefore there exists a bijection $\alpha: X\To\Gamma$. By
    the preceding observation, the components
    $\alpha_i=\pi^n_i\circ\alpha$ are onto and generous, so
    $\alpha_i\in\I_0\ss\L$ for all $1\leq i\leq n$. Since $\alpha$
    is injective, $\alpha[X]=\Gamma$ and $f\rest_\Gamma$ is
    injective, we have that $g=f(\alpha_1,\ldots,\alpha_n)\in\cl{\{f\}\cup\L}$ is
    injective. Furthermore, $g[X]=f[\Gamma]\ss\bigcup_{i=1}^n A_i$
    is co-large. Whence
    $\Oo\ss\cl{\{g\}\cup\L}\ss\cl{\{f\}\cup\L}$ by Lemma
    \ref{LEM:universalFunction} and we are done.
\end{proof}

\begin{lem}\label{LEM:columnsWithManyDifferentValuesII}
    If $f\in\On\sm\U$ is so that for all
    $1\leq i\leq n$ there exist $c\in X$ and $S\ss
    H^i_c$ such that $f[S]$ is
    large and such that for all $b\in X$ the set $\{x\in S: f(x)\neq
    f(x_1,\ldots,x_{i-1},b,x_{i+1},\ldots,x_n)\}$ is small, then $\cl{\{f\}\cup\L}\supseteq\Oo$.
\end{lem}
\begin{proof}
    Fix for every $1\leq i\leq n$ an element $c_i\in
    X$ and a set $S_i\subseteq H^i_{c_i}$
    such that $f[S_i]$ large and such that for all $b\in X$ the
    set $\{x\in S_i: f(x)\neq
    f(x_1,\ldots,x_{i-1},b,x_{i+1},\ldots,x_n)\}$ is small. Set
    $A_i=f[S_i]$, $1\leq i\leq n$. By thinning out the $S_i$ we
    can assume that the $A_i$ are disjoint and that $\bigcup_{i=1}^n A_i$ is co-large. Now one follows the proof of
    the preceding lemma.
\end{proof}

\begin{lem}
    If $f\in\On\sm\U$, then there exists $g\in\cl{\{f\}\cup\L}$ having
    co-large range and with the property that $\{x\in
    X:|g\inv[x]|=1\}$ is large (that is, the kernel of $g$ has $\kappa$ one-element classes).
\end{lem}
\begin{proof}
    There is nothing to prove if $f$ satisfies the condition of Lemma
    \ref{LEM:columnsWithManyDifferentValuesII}, so assume it does not, and let $i=1$ witness this.
    Take $c\in X$ such that
    $f[H^1_c]$ is large and choose $S\ss H^1_c$ with the property
     that $f[S]$ is still large and that $f\rest_S$ is injective.
    By assumption, there exists $b\in X$ such that $\{x\in S:
    f(x)\neq f(b,x_2,\ldots,x_n)\}$ is large. Thus, we can find a large $A\ss
    S$ with the property that $f[A]$ and $f[\{(b,x_2,\ldots,x_n): x\in
    A\}]$ are disjoint and such that the union of these two sets
    is co-large.
    Choose now generous $\alpha_2,\ldots,\alpha_n\in\Oo$ such that
    $(c,\alpha_2,\ldots,\alpha_n)[X]=A$. Since $\cl{\{f\}\cup\L}$ contains all generous functions by Lemma
    \ref{LEM:allGenerousFunctions}, we have $\alpha_j\in\L$ for $2\leq
    j\leq n$. Take a large and co-large $B\ss X$ such that
    $(c,\alpha_2,\ldots,\alpha_n)\rest_B$ is injective. Define
    $$
        \alpha_1(x)=
        \begin{cases}c&,x\in B\\
        b&,\ow
        \end{cases}
    $$
    and set $g=f(\alpha_1,\ldots,\alpha_n)$. Then $g\in\cl{\{f\}\cup\L}$ as $\alpha_1\in\X\ss\L$.
    Clearly, $(\alpha_1,\ldots,\alpha_n)\rest_B$ is injective and so is
    $g\rest_B$. Since $g[B]$ and $g[X\sm B]$ are disjoint we have
    that $|g\inv[x]|=1$ for all $x\in g[B]$. Moreover, $g[X]\ss f[A]\cup f[\{(b,x_2,\ldots,x_n): x\in
    A\}]$ is co-large.
\end{proof}

\begin{lem}
    Let $f\in\On\sm \U$. If $h\in\Oo$ is a function whose
    kernel has at least one large equivalence class (that is, there
    exists $x\in X$ with $h\inv[x]$ large), then $h\in \cl{\{f\}\cup\L}$.
\end{lem}
\begin{proof}
    There exist a large $B\ss X$ and $b\in X$ such that
    $h[B]=\{b\}$. Let $g$ be provided by the preceding lemma.
    With the help of permutations of the base set we can assume that $|g\inv
    [x]|=1$ for all $x\in g[X\sm B]$. Since the range of $g$
    is co-large we can find $\delta: X\sm g[X] \To X$ onto and
    generous. Now define $m\in\Oo$ by
    $$
        m(x)=
        \begin{cases}
            \delta(x)&,x\nin g[X]\\
            b&, x\in g[B]\\
            h(g\inv(x))&, x\in g[X\sm B].
        \end{cases}
    $$
    Obviously $m\in \I_0\ss\L$ and $h=m\circ g\in\cl{\{f\}\cup\L}$.
\end{proof}

Having found many functions which $\cl{\{f\}\cup\L}$ must contain,
we are finally ready to prove Lemma \ref{LEM:LandNotAlmostUnary}.

\begin{proof}[Proof of Lemma \ref{LEM:LandNotAlmostUnary}]
    There are $c_1,\ldots,c_n\in X$ such
    that $f[H^i_{c_i}]$ is large for $1\leq
    i\leq n$. Take for all $1\leq i\leq n$ large $B_i\subseteq H^i_{c_i}$ with the property
    that $f\rest_B$ is injective and $f[B]$ is co-large,
    where $B={\bigcup_{i=1}^n
    B_i}$. Let $\alpha: X\To B$ be any bijection. Since $\alpha_i\inv[c_i]$ is large for every
    component $\alpha_i=\pi^n_i\circ\alpha$, the preceding lemma yields $\alpha_i\in\cl{\{f\}\cup\L}$ for
    $1\leq i\leq n$. Whence, $g=f(\alpha_1,\ldots,\alpha_n)\in\cl{\{f\}\cup\L}$. But
    $g[X]=f[B]$ is co-large and $g$ is injective by construction;
    thus Lemma \ref{LEM:universalFunction} yields $\Oo\ss\cl{\{g\}\cup\L}\ss\cl{\{f\}\cup\L}$.
\end{proof}

\begin{proof}[Proof of Proposition \ref{PROP:generalizedGavrilov}]
    Since $\cl{\Oo\cup\H}=\O$, there must exist some
    $f\in\H\sm\U$. But then, since $\G\supseteq \L$, Lemma \ref{LEM:LandNotAlmostUnary} implies
    $\cl{\G\cup\H}\supseteq\Oo$ so that we
    infer $\cl{\G\cup\H}=\O$.
\end{proof}

\end{section}

\begin{section}{The proof of Theorem \ref{THM:allMaximalMonoids}}

We now determine on an infinite $X$ all maximal submonoids of
$\Oo$ which contain the permutations, proving Theorem
\ref{THM:allMaximalMonoids}. In a first section, we present the
part of the proof which works on all infinite sets; then follow
one section specifically for the case of a base set of regular
cardinality and another section for the singular case. Throughout
all parts we will mention explicitly whenever a statement is true
only on $X$ of regular or singular cardinality, respectively.

\subsection{The part which works for all infinite sets}
\begin{prop}\label{PROP:GlambdaIsMaximal}
    $\G_\lambda$ is a maximal submonoid of $\Oo$ for $\lambda=1$ and
    $\aleph_0\leq\lambda\leq\kappa$.
\end{prop}
\begin{proof}
    As already mentioned in Lemma \ref{LEM:GlambdaFacts}, the
    maximality of the $\G_\lambda$ for $\lambda=1$ or infinite has
    been proved in \cite{Ros74} (Lemma 5.2).
\end{proof}

The maximal monoids of Proposition \ref{PROP:GlambdaIsMaximal}
already appeared in the preceding section since they give rise to
maximal clones via $\pol$. We shall now expose maximal monoids
above the permutations which do not have this property. Recall
that $\M_\lambda$ consists of all functions which are either
$\lambda$-surjective or not $\lambda$-injective.

\begin{prop}\label{PROP:MlambdaIsMaximal}
    Let $\lambda=1$ or $\aleph_0\leq\lambda\leq\kappa$. Then $\M_\lambda$ is a maximal submonoid of $\Oo$.
\end{prop}
\begin{proof}
    We show first that $\M_\lambda$ is closed under composition.
    Let therefore $f,g\in\M_\lambda$, that is, those functions are
    either $\lambda$-surjective or not $\lambda$-injective; we claim that $f\circ g$ has either of these properties. It is
    clear that if $g$ is not $\lambda$-injective, then $f\circ g$
    has the same property. So let $g$ be $\lambda$-surjective. It is easy to see that if
    $f$ is $\lambda$-surjective, then so is $f\circ g$. So
    assume finally that $f$ is not $\lambda$-injective. We claim
    that $f\circ g$ is not $\lambda$-injective either. For
    $\lambda=1$ this is just the statement that if $f$ is not
    injective, and $g$ is surjective, then $f\circ g$ is not
    injective, which is obvious. Now consider the infinite case.
    There exist disjoint $A,B\ss X$ of size $\lambda$ such that
    $f[A]=f[B]$. Set $A'=A\cap g[X]$; $A'$ still has size
    $\lambda$ as $g$ misses less than $\lambda$ values. Clearly
    $B'=\{x\in B:\exists y\in A' (f(x)=f(y))\}$ has size $\lambda$
    as well and so does $B''=B'\cap g[X]$. But now for the sets
    $C=g\inv [A']$ and $D=g\inv [B'']$ it is true that $|C|,|D|
    \geq\lambda$, $C\cap D=\emptyset$, and $f\circ g[C]=f\circ
    g[D]$; hence $f\circ g$ is not $\lambda$-injective.\\
    Now we prove that $\M_\lambda$ is maximal in $\Oo$. Consider for
    this reason any $m\nin\M_\lambda$, that is, $m$ is
    $\lambda$-injective and misses at least $\lambda$ values. There
    exists $A\ss X$ so that $|X\sm A|<\lambda$ and such that the restriction of $m$ to $A$ is injective.
    Take any injection $i\in\Oo$ with $i[X]=A$. Then $i\in\M_\lambda$ as $i$ is $\lambda$-surjective. Now let
    $f\in\Oo$ be arbitrary. Define
    $$
        g(x)=
        \begin{cases}
            f((m\circ i)\inv(x))&,x\in m\circ i[X]\\
            a&, \ow\\
        \end{cases}
    $$
    where $a\in X$ is any fixed element of $X$. Being constant on the
    complement of the range of $m$, $g$ it is not $\lambda$-injective
    and whence an element of $\M_\lambda$. Therefore $f=g\circ
    m\circ i\in\cl{\M_\lambda\cup\{m\}}$ so that we infer
    $\cl{\M_\lambda\cup\{m\}}\supseteq\Oo$.
\end{proof}

\begin{lem}\label{LEM:allMonoidsAboveI0}
    There are no other maximal monoids above $\S\cup\I_0$ except the
    $\M_\lambda$ ($\lambda=1$ or $\aleph_0\leq \lambda\leq \kappa$).
\end{lem}

\begin{proof}
    Let $\G\supseteq \I_0\cup\S$ be a submonoid of $\Oo$ which is not contained in any of the
    $\M_\lambda$; we prove that $\G=\Oo$. To do this, we
    show that $\G$ contains an injective function $u\in
    \Oo$ with co-large range; then the lemma follows from Lemma
    \ref{LEM:universalFunction}. Fix for every $\lambda$ a
    function $m_\lambda\in\G\sm\M_\lambda$. Since $m_{\kappa}$ is
    $\kappa$-injective, there exists a cardinal $\lambda_1
    <\kappa$ and a set $A_1\ss X$ of size $\lambda_1$ such that
    the restriction of $m_\kappa$ to the complement of $A_1$ is
    injective. If $\lambda_1$ is infinite, then consider
    $m_{\lambda_1}$. Not being an element of $\M_{\lambda_1}$, $m_{\lambda_1}$ misses at least
    $\lambda_1$ values. Hence by adjusting it with a suitable
    permutation we can assume that $m_{\lambda_1}[X]\ss X\sm A_1$.
    There exists a cardinal $\lambda_2<\lambda_1$ and a subset
    $A_2$ of $X$ of size $\lambda_2$ such that the restriction of
    $m_{\lambda_1}$ to the complement of $A_2$ is injective.
    Hence, writing $\lambda_0=\kappa$ we obtain that
    $m_{\lambda_0}\circ m_{\lambda_1}\in\G$ is injective on $X\sm A_2$ and misses $\kappa$
    values. We can iterate this to arrive after a finite number of steps at a set
    $A_n$ of finite size $\lambda_n$ such that the restriction of
    $m_{\lambda_{0}}\circ \ldots\circ m_{\lambda_{n-1}}\in\G$ to $X\sm
    A_n$ is injective and misses $\kappa$ values. Since
    $m_1\nin\M_1$ is injective and misses at least one value we
    conclude that the iterate $m_1^{\lambda_n}\in\G$ is injective and
    misses at least $\lambda_n$ values. Modulo
    permutations we may assume that $m_1^{\lambda_n}[X]\ss X\sm
    A_n$. But now we have that $m_{\lambda_{0}}\circ \ldots\circ
    m_{\lambda_{n-1}}\circ m_1^{\lambda_n}\in \G$ is injective and misses $\kappa$ values,
    implying that $\G=\Oo$.
\end{proof}

\subsection{The case of a base set of regular cardinality}

We now finish the proof of Theorem \ref{THM:allMaximalMonoids} for
the case when $X$ has regular cardinality. The proof for this case
comprises Propositions \ref{PROP:GlambdaIsMaximal},
\ref{PROP:MlambdaIsMaximal}, \ref{PROP:AisMaximal} and
\ref{PROP:noOtherMaxMonoids}.

\begin{prop}\label{PROP:AisMaximal}
    If $X$ is of regular cardinality, then $\A$ is a maximal submonoid of $\Oo$.
\end{prop}

\begin{proof}
    This has been proved in \cite{Ros74} (Proposition 4.1).
\end{proof}

\begin{prop}\label{PROP:noOtherMaxMonoids}
    Let $X$ have regular cardinality. There exist no other maximal submonoids of $\Oo$ containing the permutations
    except those listed in Theorem \ref{THM:allMaximalMonoids} for the regular case.
\end{prop}
\begin{proof}
    Assume that $\G\supseteq \S$ is a submonoid of $\Oo$ not contained in any
    of the monoids of the theorem; we show that $\G=\Oo$. Indeed, since
    $\G\nsubseteq\A$, Proposition \ref{PROP:someIlambdaIsBelow} tells us that there exists a cardinal $\lambda
    \leq\kappa$ such that $\I_\lambda$ is contained in $\G$. Choose
    $\lambda$ minimal with this property. If $\lambda$ was greater than $0$, then
    $\G\ss\G_\lambda$ for otherwise Lemma \ref{LEM:existsLambda0}
    would yield a contradiction to the minimality of $\lambda$. But
    this is impossible as we assumed that $\G$ is not contained in
    any of the $\G_\lambda$, so we conclude that $\lambda=0$. Now
    Lemma \ref{LEM:allMonoidsAboveI0} implies that $\G=\Oo$.
\end{proof}

\subsection{The case of a base set of singular cardinality}

The only problem with base sets of singular cardinality is that
the set $\A$ is not closed under composition; in fact,
$\cl{\A}=\O$. A slight adjustment of the definition of $\A$ works
in this case. We will refer to results from preceding sections;
this might look unsafe since there we restricted ourselves to base
sets of regular cardinality. However, when proving the particular
results cited here we did not use the regularity of the base set.
The proof of Theorem \ref{THM:allMaximalMonoids} for singular
cardinals comprises Propositions \ref{PROP:GlambdaIsMaximal},
\ref{PROP:MlambdaIsMaximal}, \ref{PROP:A'isMaximal} and
\ref{PROP:noOtherMaximalMonoidsSingular}.

\begin{defn}
    A function $f\in\Oo$ is said to be \emph{harmless} iff
    there exists $\lambda < \kappa$ such that the set of all $x\in X$ for which
    $|f\inv [x]|>\lambda$ is small. With this definition,
    $\A'$ as defined in Theorem \ref{THM:allMaximalMonoids} is the set of all harmless functions.
\end{defn}

\begin{lem}
    $\A'$ is a monoid and $\A'\ss \A$. Moreover, $\A=\A'$ iff $\kappa$ is a
    successor cardinal.
\end{lem}
\begin{proof}
    It is obvious that $\A'\ss \A$ and that $\A=\A'$ iff $\kappa$ is a
    successor cardinal. To prove that $\A'$ is closed
    under composition, let $f,g\in\A'$; we show $h=f\circ g\in \A'$.
    There exist $\lambda_f,\lambda_g<\kappa$ witnessing that $f$
    and $g$ are harmless. Set $\lambda$ to be $\max(\lambda_f,\lambda_g)$; we
    claim that the set of $x\in X$ for which $|h\inv [x]|>\lambda$
    is small. For if $|h\inv [x]|>\lambda$, then either $|g\inv
    [x]|>\lambda$ or there exists $y\in g\inv
    [x]$ such that $|f\inv [y]|>\lambda$. Both possibilities occur
    only for a small number of $x\in X$ and so $h$ is harmless.
\end{proof}

\begin{lem}\label{LEM:notInA'MakesNotInA}
    Let $X$ have singular cardinality. If $g\nin\A'$, then $g$ together with $\S$ generate a function not in $\A$.
\end{lem}

\begin{proof}
    Set $\lambda < \kappa$ to be the cofinality of $\kappa$. Because $g$ is not harmless, there
    exist distinct sequences $(x^0_\xi)_{\xi < \lambda},\ldots,(x^\kappa_\xi)_{\xi < \lambda}$ of distinct elements of $X$
    such that $\bigcup_{\xi <\lambda} g\inv [x^\zeta_\xi]$
    is large for all $\zeta <\kappa$. Indeed, if $(\mu_\xi)_{\xi<\lambda}$ is any cofinal sequence of cardinalities
    in $\kappa$, then the fact that $g$ is not harmless allows us to pick for every $\xi<\lambda$ an element $x^0_\xi\in X$ such
    that $|g\inv[x^0_\xi]|>\mu_\xi$; it is also no problem to choose the elements distinct. This yields the first sequence and
    since with every sequence we are using up only $\lambda <\kappa$ elements, the definition of harmlessness ensures that we can
    repeat the process $\kappa$ times. By throwing away half of the sequences, we may assume that the set of all $y\in X$ which
    do not appear in any of the sequences is large.\\
    There exists a permutation $\alpha\in\S$ such that
    $g\circ\alpha(x_{\xi_1}^{\zeta_1})=g\circ\alpha (x_{\xi_2}^{\zeta_2})$ if and only if $\zeta_1=\zeta_2$,
    for all $\zeta_1,\zeta_2<\kappa$ and all $\xi_1,\xi_2
    <\lambda$. For we can map every sequence $(x_\xi^\zeta)_{\xi <\lambda}$ injectively into an equivalence
    class of the kernel of $g$ of size greater than $\lambda$;
    since there are many such classes every sequence can be assigned an
    own class, and we choose the classes so that a large number of classes are not hit at all.
    This partial injective mapping we can then extend to the permutation $\alpha$ as
    it is defined on a co-large set and has co-large range.\\
    Set $y^\zeta=g\circ\alpha(x_0^\zeta)$ for all $\zeta<\kappa$. Then
    the
    $y^\zeta$ are pairwise distinct and for all $\zeta<\kappa$ we have that
    $(g\circ\alpha\circ g)\inv[y^\zeta]\supseteq \bigcup_{\xi <\lambda} g\inv [x^\zeta_\xi]$
    is large. Hence, $g\circ \alpha \circ g\nin
    \A$.
\end{proof}

\begin{prop}\label{PROP:A'isMaximal}
    Let $X$ have singular cardinality. Then $\A'$ is a maximal submonoid of $\Oo$.
\end{prop}

\begin{proof}
    Let $g\in\Oo\sm\A'$. We know that $g$ together with $\A'$ generate a function not in $\A$.
    Then by Lemma \ref{LEM:generousInfiniteRange}, we obtain a function which is generous and has large range,
    call it $h$. Now take any $f\in\Oo$ such that $f\circ h[X]=X$
    which is injective on $h[X]$ and constant on $X\sm h[X]$. Then
    $f\in\A'$ and $f\circ h\in\I_0$. Thus,
    $\I_0\ss\cl{\{g\}\cup\A'}$ and since all injections are elements of $\A'$ we can apply
    Lemma \ref{LEM:universalFunction} to prove
    $\cl{\{g\}\cup\A'}\supseteq\Oo$.
\end{proof}

\begin{prop}\label{PROP:noOtherMaximalMonoidsSingular}
    Let $X$ have singular cardinality. There
    exist no other maximal submonoids of $\Oo$ containing the permutations
    except those listed in Theorem \ref{THM:allMaximalMonoids} for the singular case.
\end{prop}

\begin{proof}
    If $\G\supseteq\S$ is a submonoid of $\Oo$ which is not
    contained in $\A'$, then it is not contained in $\A$ by Lemma \ref{LEM:notInA'MakesNotInA}. From
    this point, one can follow
    the proof of Proposition \ref{PROP:noOtherMaxMonoids}.
\end{proof}

\end{section}

\providecommand{\bysame}{\leavevmode\hbox
to3em{\hrulefill}\thinspace}
\providecommand{\MR}{\relax\ifhmode\unskip\space\fi MR }
\providecommand{\MRhref}[2]{%
  \href{http://www.ams.org/mathscinet-getitem?mr=#1}{#2}
} \providecommand{\href}[2]{#2}

\end{document}